\documentclass[12pt]{elsart}

\usepackage{lscape}
\usepackage{longtable,booktabs}
\usepackage{amssymb}
\usepackage{amsmath}
\begin{document}

\title{Complete Group Classification of Systems of Two
\\
Nonlinear Second-Order
Ordinary Differential
\\ Equations of the Form ${\bf y}^{\prime\prime}={\bf F}({\bf y})$}

\author[SUT]{G.F. Oguis}
\ead{fae@math.sut.ac.th}

\author[DUT]{S. Moyo}
\ead{moyos@dut.ac.za}

\author[SUT]{S.V. Meleshko}
\ead{sergey@math.sut.ac.th}

\address[SUT]{Suranaree University of Technology, School of Mathematics, Nakhon
Ratchasima 30000, Thailand}

\address[DUT]{Durban University of Technology, Research and Postgraduate Support \& Institute for Systems Science, P O Box 1334, Steve
Biko Campus, Durban 4000, South Africa}


\begin{frontmatter}



\begin{abstract}
Extensive work has been done on the group classification of systems of equations in the literature. This paper identifies the gap in the literature which concerns the group classification of systems of two nonlinear second-order ordinary differential equations. We provide a complete group classification of systems of two ordinary differential equations of the form, $\bf{y}^{\prime\prime}=\bf{F}(\bf{y})$, which occur in many physical applications using two approaches which form the essence of this paper.
\end{abstract}

\begin{keyword}
Group classification\sep autonomous systems\sep nonlinear systems\sep admitted Lie group
\end{keyword}

\end{frontmatter}

\section{Introduction}

Systems of autonomous nonlinear second-order ordinary differential equations, where the independent variable, usually assumed to be time, does not appear on the right hand side of the system, arise in various physical problems. This effectively assumes that the laws of nature which hold true in the present are presumably applicable in the past and the future. Hence, like all other systems of differential equations, the study of their symmetry structure poses an important role as their presence in a system allows one to reduce the order of the studied equations and also to find general solutions in quadratures.

Group classification studies, dating more than a century back, were first initiated by the founder of symmetry analysis, Sophus Lie \cite{bk:Lie[1883],bk:Lie[1891b],bk:Lie[1884],bk:Lie[1881]}. These studies were long forgotten until Ovsiannikov \cite{bk:ovsiannikov[1958],bk:Ovsiannikov[1978]} revived the work five decades ago. Lie's works put emphasis on tackling the group classification in two ways: 1) the direct way and 2) the indirect way also known as the algebraic approach. The direct way involves directly finding solutions of the determining equations and allows one to study all possible admitted Lie algebras without omission. On the other hand the indirect way involves solving the determining equations up to finding relations between constants defining admitted generators. The algebraic approach, as in the studies \cite{bk:MahomedLeach[1989],bk:GonzalezKamranOlver[1992a],bk:PopovychKunzingerEshraghi[2010],bk:GrigorievMeleshkoSuriyawichitseranee[2013]}\footnote{See also references therein.}, takes into account the algebraic properties of an admitted Lie group and the knowledge of the algebraic structure of the admitted Lie algebras in order to allow significant simplification of the group classification.
In one of Lie's works \cite{bk:Lie[1883]}, he gave a complete group classification of a single second-order ordinary differential equation of the form $y^{\prime\prime}=f(x,y)$. Later on  Ovsiannikov \cite{bk:Ovsiannikov[2004]} did this group classification in a different way. The method he used, now also known as the direct approach, involved a two-step technique where the determining equations were first simplified through exploiting equivalence transformations and later on solved for the reduced cases of the generators. The same technique was used in a study conducted in \cite{bk:Phauk[2013]} to classify a more general case of equations of the form $y^{\prime\prime}=P_3(x,y;y^\prime),$ where $P_3(x,y;y^\prime)$ is a polynomial of a third degree with respect to the first-order derivative $y^\prime.$ Observe that sometimes difficulties arise in using the direct approach. Sometimes it is difficult to select or tease out equivalent cases with respect to equivalence transformations. As is observed in the classification of a general scalar second-order ordinary differential equation of the form $y^{\prime\prime}=f(x,y;y^\prime),$ the application of the direct technique gives rise to overwhelming difficulties. In this study both the direct and indirect techniques are employed.

Apart from dealing with classification problems there is a significant amount of research that deals with the dimension and structure of symmetry algebras of linearizable ordinary differential equations \cite{bk:MahomedLeach[1989],bk:MahomedLeach[1990],bk:HandbookLie_v3,bk:Gorringeleach[1988],bk:WafoMahomed[2000],bk:BoyPopoSha[2012]}. This is also of importance since some nonlinear equations appear in disguised forms.

In addition to extensive studies on properties of scalar second-order ordinary differential equations, there are also several researchers committed to studying systems of two linear second-order ordinary differential equations \cite{bk:WafoSoh[2010],bk:Meleshko[2011],bk:BoyPopoSha[2012],%
bk:Campoamor-Stursberg[2011],bk:Campoamor-Stursberg[2012],%
bk:MoyoMeleshkoOguis[2013],bk:MkhizeMoyoMeleshko[2014],%
bk:MoyoMeleshkoOguis[2014],bk:MeleshkoMoyo[2014]}.
 Surprisingly, the group classification of systems of two nonlinear second-order ordinary differential equations has not yet been exhausted, in particular, the group classification of systems of two autonomous nonlinear second-order ordinary differential equations is not yet complete. Hence this paper considers the group classification of systems of two autonomous nonlinear second-order ordinary differential equations of the form
\begin{equation}\label{eq:auto}\bf{y}^{\prime\prime}=\bf{F}(\bf{y}).
\end{equation}

The system studied here is a generalization of Lie's study \cite{bk:Lie[1891b]}. Studied cases such as systems of two linear second-order ordinary differential equations and the degenerate case which is equivalent to the following
\begin{equation}\label{eq:dege}
y^{\prime\prime}=F(x,y,z),\quad\quad z^{\prime\prime}=0
\end{equation} are omitted from this paper. We call systems that are equivalent to these cases as reducible systems and irreducible otherwise.

The paper is organized as follows. A preliminary study of systems of two nonlinear second-order ordinary differential equations is tackled first and is followed by the subsequent group classification applied to autonomous systems (\ref{eq:auto}) of two second-order ordinary differential equations. The group classification is divided into two parts depending on the coefficient of the infinitesimal generator. The direct approach is applied on one case while a combination of the optimal system of subalgebras and direct approach is applied to the other case. The latter part of the paper lists the different cases with their respective results and is then followed by the conclusion.

\section{Background study of systems of the form $\mathbf{y}^{\prime\prime}=\mathbf{F}(x,\mathbf{y})$}

This section focuses on systems of two nonlinear second-order ordinary differential equations of the form \cite{bk:MoyoMeleshkoOguis[2013],bk:MeleshkoMoyo[2014]}
\begin{equation}\label{eq:1}
\mathbf{y}^{\prime\prime}=\mathbf{F}(x,\mathbf{y}),
\end{equation}
where
\[
\mathbf{y} = \left(
\begin{array}{c}
y  \\
z %
\end{array}
\right),\quad\quad \mathbf{F}=
\left(
\begin{array}{c}
F(x,y,z)  \\
G(x,y,z) %
\end{array}
\right).
\]

\subsection{Equivalence transformations}
System (\ref{eq:1}) has the following equivalence transformations:
\begin{enumerate}
\item \label{linchan1} a linear change of the dependent variables $\mathbf{\tilde{y}}=P\mathbf{y}$ with constant nonsingular $2\times 2$ matrix $P$;
\item the change
$\tilde{y}=y+\phi(x)$ and $\tilde{z}=z+\psi(x);$ and
\item the transformation related with the change $\tilde{x}=\phi(x),$ $\tilde{y}=y\psi(x),$ $\tilde{z}=z\psi(x),$
where the functions $\phi(x)$ and $\psi(x)$ satisfy the condition
$\displaystyle \frac{\phi^{\prime\prime}}{\phi^{\prime}} = 2 \frac{\psi^{\prime}}{\psi}. $
\end{enumerate}
%

\subsection{Determining equations}
The determining equations in matrix form for irreducible systems of the form (\ref{eq:1}) are given by
\begin{equation}\label{eq:2}
2\xi\mathbf{F}_x + 3\xi^\prime\mathbf{F} +(((A+\xi^\prime E)\mathbf{y} + \zeta)\cdot\nabla)\mathbf{F} - A\mathbf{F} = \xi^{\prime\prime\prime}\mathbf{y} + \zeta^{\prime\prime},
\end{equation}
where the matrix $A = (a_{ij})$ is constant. The associated infinitesimal generator is
$$
\begin{array}{c}
X = 2 \xi(x) \partial_x + (A\mathbf{y} +\zeta(x))\cdot\nabla,
\end{array}
$$
where $\nabla=(\partial_y,\partial_z)^t$ and $"\cdot"$ means the scalar product $\mathbf{b}\cdot\nabla = b_i \partial_{y_i},$ where the summation with respect to the repeated index is used \cite{bk:MoyoMeleshkoOguis[2013]}.

The equivalence transformation (\ref{linchan1}) with linear change $\mathbf{\tilde{y}}=P\mathbf{y}$, when applied to equations (\ref{eq:1}), reduces equations (\ref{eq:2}) and its associated infinitesimal generator to the same form with the matrix $A$ and the vector $\zeta$ changed. Equations (\ref{eq:1}) become
$$
\mathbf{\tilde{y}}^{\prime\prime} = \mathbf{\tilde{F}}(x,\mathbf{\tilde{y}})
$$ with
$$
\mathbf{\tilde{F}}(x,\mathbf{\tilde{y}})  = P \mathbf{F}(x,P\mathbf{\tilde{y}}),
$$
and the partial derivatives with respect to the variables $\mathbf{y}$ are also changed as follows
$$
\mathbf{b}\cdot \nabla = (P \mathbf{b}) \cdot \tilde{\nabla}.
$$
Consequently, the determining equations (\ref{eq:2}) become
\[
2\xi\mathbf{\tilde{F}}_x + 3\xi^\prime\mathbf{\tilde{F}} +(((\tilde{A}+\xi^\prime E)\mathbf{\tilde{y}} + \tilde{\zeta})\cdot\tilde{\nabla})\mathbf{\tilde{F}} - \tilde{A}\mathbf{\tilde{F}} - \xi^{\prime\prime\prime}\mathbf{\tilde{y}} - \tilde{\zeta}^{\prime\prime} = 0,
\]
where
$$
\tilde{A} = P A P^{-1},\quad\quad \tilde{\zeta} = P \zeta
$$
and the associated infinitesimal generator is also changed as follows
$$
\begin{array}{c}
X = 2 \xi(x) \partial_x + (\tilde{A}\tilde{\mathbf{y}} +\tilde{\zeta(x)})\cdot\tilde{\nabla}.
\end{array}
$$
 As in \cite{bk:MoyoMeleshkoOguis[2013],bk:MeleshkoMoyo[2014]} and in the succeeding pages, this transformation places a very important role in the group classification process.

From the study \cite{bk:MoyoMeleshkoOguis[2013]}, the systems of two nonlinear second-order ordinary differential equations are equivalent to one of the following 10 types listed below. Looking closely at these systems, there is a necessity to conduct an initial study where the systems of two equations do not depend on $x.$ This forms the core of this paper.
\pagebreak
\begin{center}
\begin{table}[!htb]
\begin{tabular}{cll}
   & \underline{$F$ and $G$} &\underline{Relations and conditions} \\
 1. & $F=e^{ax}f(u,v)$,  & $u=ye^{-ax},\,\,v=ze^{-bx},$ \\
   & $G=e^{bx}g(u,v)$   & $a,\,\,b$ are constant \\
 2. & $F= e^{ax}(\cos(cx)f(u,v)+\sin(cx)g(u,v)),$ & $u=e^{-ax}(y \cos(cx)-z\sin(cx)),$ \\
   & $G=e^{ax}(-\sin(cx)f(u,v)+\cos(cx)g(u,v))$ & $v=e^{-ax}(y\sin(cx)+z\cos(cx)),$\\
   & & $a,\,c\neq0$  are constant\\
 3. & $F= e^{ax}(f(u,v)+xg(u,v)),$ & $u=e^{-ax}(y-zx),\,\,v=ze^{-ax},$ \\
   & $G=e^{ax}g(u,v)$ & $a$ is constant\\
 4. & $F= (y+h_1(x))f(x,v)-h_1^{\prime\prime},$ &$v=(z+h_2(x))(y+h_1(x))^{\alpha},$ \\
   & $G=(z+h_2(x))g(x,v)-h_2^{\prime\prime}$ & $\alpha\neq0$ is constant\\
 5. & $F=(y+h_1(x))f(x,v)-h_1^{\prime\prime}(x),$ & $v=z- h_2(x)\ln{(y+h_1(x))}$\\
   & $G=h_2^{\prime\prime}(x)\ln{(y+h_1(x))}+g(x,v)$ & \\
 6. & $F= \displaystyle{\frac{h_1^{\prime\prime}(x)}{h_1(x)}}y+f(x,v),$ & $v=z- \displaystyle{\frac{h_2(x)}{h_1(x)}}y,\,\,h_1(x)\neq0$ \\
   & $G= \displaystyle{\frac{h_2^{\prime\prime}(x)}{h_1(x)}}y+g(x,v)$ & \\
 7. & $F= e^{au}(\cos(cu)f(x,v)+\sin(cu)g(x,v)),$ & $y=ve^{au}\sin(cu),\,\,z=e^{au}\cos(cu),$\\
   & $G=e^{au}(-\sin(cu)f(x,v)+\cos(cu)g(x,v))$ & $a,\,c\neq0$ are constant\\
 8. & $F = \displaystyle{\frac{y}{z+h_1(x)}}f(x,v)+g(x,v),$ & $v=z+ h_1(x)$\\
   & $G=-h_1^{\prime\prime}(x)+f(x,v)$ & \\
 9. & $F= \displaystyle{\frac{h_2^{\prime\prime}(x)}{2}}u^2+uf(x,v)+g(x,v),$ & $u= \displaystyle{\frac{z+h_1(x)}{h_2(x)}},\,\,v=y-\displaystyle{\frac{(z+h_1(x))^2}{h_2(x)}},$ \\
   & $G=-h_1^{\prime\prime}(x)+h_2^{\prime\prime}(x)u+f(x,v)$ & $h_2(x)\neq0$\\
 10.& $F= e^{u}(uf(x,v)+g(x,v)),$ & $y=uve^{u},\,\,z=ve^{u}$ \\
   & $G=e^{u}f(x,v)$ & \\
\end{tabular}
\begin{tabular}{l}
For all the cases, $h_1$, $h_2$, $f$ and $g$ are arbitrary functions of their arguments.\\
\end{tabular}
\end{table}
\end{center}

\section{Autonomous systems (\ref{eq:auto}) of two nonlinear second-order ordinary differential equations and their group classification}
Since for autonomous systems $\mathbf{F}_{x}=0$, the determining equations for autonomous systems have the form
\begin{equation}\label{eq:detanl}
3\xi^\prime\mathbf{F} +(((A+\xi^\prime E)\mathbf{y} + \zeta)\cdot\nabla)\mathbf{F} - A\mathbf{F} - \xi^{\prime\prime\prime}\mathbf{y} - \zeta^{\prime\prime} = 0.
\end{equation} This also implies that the generator $\partial_x$ is admitted by system (\ref{eq:auto}).

Differentiating the determining equations (\ref{eq:detanl}) with respect to $x$, the group classification study is reduced into two cases, namely,
\begin{enumerate}
\item the case with at least one admitted generator with $\xi^{\prime\prime}\neq0$; and
\item the case where all admitted generators have $\xi^{\prime\prime}=0$.
\end{enumerate}
For the first case, the direct approach by Lie is utilized, whereas for the second case, a combination of the optimal system of subalgebras of the Lie algebra and the direct method is used.

\subsection{Systems admitting at least one admitted generator with $\xi^{\prime\prime}\neq0$}
For the case with at least one generator, $\xi^{\prime\prime}\neq0$,  we initially consider the differentiated determining equations (\ref{eq:detanl}) with respect to $x$ and divide them by $\xi^{\prime\prime}.$ The determining equations become
\begin{equation}\label{eq:detx}
\begin{array}{c}
\displaystyle{3\mathbf{F} + \left(\left(\mathbf{y} +\frac{\zeta^\prime}{\xi^{\prime\prime}}\right)\cdot \nabla\right)\mathbf{F} - \frac{\xi^{(4)}}{\xi^{\prime\prime}}\mathbf{y} - \frac{\zeta^{\prime\prime\prime}}{\xi^{\prime\prime}} = 0}.
\end{array}
\end{equation}
Fixing $x$, and shifting $y$ and $z$, equations (\ref{eq:detx}) are reduced to
\[
\begin{array}{c}
3\mathbf{F} + (\mathbf{y}\cdot \nabla)\mathbf{F} - a \mathbf{y} - \mathbf{b} = 0
\end{array}
\]
where vector $\mathbf{b} = (b,c)^t,$ with $a, b, c$ constant.

The general solution of these equations is
\begin{equation} \label{eq:auto2}
\begin{array}{c}
F = \displaystyle{\frac{b}{3} + \frac{a y}{4} + y^{-3} f(u)},\\
G = \displaystyle{\frac{c}{3} + \frac{a z}{4} + z^{-3} g(u)},\\ %
\end{array}
\end{equation}
where $u=z/y$ and $f^{\prime}g^{\prime} \neq 0.$ It is easy to see that if $f^{\prime}g^{\prime}=0$, the studied system is equivalent to a reducible case. The functions $F$ and $G$ are then substituted into the determining equations (\ref{eq:detanl}). The determining equations are then solved directly in order to find generators admitted by equations (\ref{eq:auto}). 
The first part of the determining equations is given as follows:
\begin{equation}\label{eq:detxn}
\begin{array}{l}
  \xi^{\prime\prime\prime} - a \xi^{\prime} =0,\\
\end{array}
\end{equation}
\begin{equation}\label{eq:detxn2}
\begin{array}{c}
(\zeta_1 u - \zeta_2) f_u + 3 \zeta_1 f =0,\\
(u^2 \zeta_1 -  u \zeta_2)g_u + 3 \zeta_2 g =0,\\
\end{array}
\end{equation}
\begin{equation}\label{eq:detxn3}
\begin{array}{c}
 12 \zeta_1^{\prime\prime}- 12 b \xi^{\prime} - 3 a \zeta_1 + 4 a_{11} b + 4 a_{12} c=0,\\
 12 \zeta_2^{\prime\prime}- 12 c \xi^{\prime}  - 3 a \zeta_2 + 4 a_{21} b + 4 a_{22} c=0,\\
\end{array}
\end{equation}
\begin{equation}\label{eq:detxn4}
\begin{array}{c}
 (a_{11} u^4 + a_{12} u^5 - a_{21} u^3 - a_{22} u^4) f_u
+ (4 a_{11}  u^3 + 3 a_{12} u^4) f + a_{12} g=0,\\
 (a_{11} u^2 +  a_{12} u^3 - a_{21} u - a_{22} u^2) g_u
+ a_{21} u^4 f + (3 a_{21} + 4 a_{22} u) g = 0.\\
\end{array}
\end{equation}

From equation (\ref{eq:detxn}), it can be seen that the general solution of $\xi$ depends on three values of $a,$ i.e., $a=0,$ $a = -p^2$ and $a = p^2$, where $p\neq0.$ 
For $a=0,$ the general solution of $\xi$ is \[\xi=\xi_2 x^2 + \xi_1 x + \xi_0,\] where $\xi_2\neq0,\,\xi_1,\,\xi_0$ are constant. For $a = -p^2,$ the general solution of $\xi$ is \[\xi=\xi_1 \cos(px) +\xi_2 \sin(px) +\xi_0,\] where $\xi_2\neq0,\,\xi_1\neq0,\,\xi_0$ are constant. Lastly, for $a = p^2,$ the general solution of $\xi$ is \[\xi=\xi_1 e^{-px} +\xi_2 e^{px} +\xi_0,\] where $\xi_2\neq0,\,\xi_1\neq0,\,\xi_0$ are constant.
Subsequently the determining equations (\ref{eq:detxn2}) lead to the study of two cases where: (1) there exists a generator with $\zeta_1\neq0$ and (2) all generators have $\zeta_1=0.$

 Considering the case where there exists a generator with $\zeta_1\neq0$, we divide by $\zeta_1$ and differentiate the equations (\ref{eq:detxn2}) with respect to $x$ to obtain, $\zeta_2 = k \zeta_1$, where $k$ is a constant.  Substituting this back to equations (\ref{eq:detxn2}), one obtains $f=f_0(u-k)^{-3}$ and $g=g_0u^3(u-k)^{-3}$. Also, differentiating equations (\ref{eq:detxn3}) with respect to $x$, it follows that $c=kb.$ From here, one can verify that this is a reducible case. 

 Consider that all generators have $\zeta_1=0.$ From equations (\ref{eq:detxn2}), it follows that $\zeta_2=0.$ Differentiating equations (\ref{eq:detxn3}) with respect to $x$, it immediately follows that $b=c=0$. From equations (\ref{eq:detxn4}), the equivalence transformation $\tilde{y}=Py,$ where $P$ is a constant nonsingular $2\times2$ matrix, is utilized in order to obtain the general solution of $f$ and $g$. Note that the constant matrix $A$ is reduced to one of the following real-valued Jordan forms
\begin{equation}\label{eq:jor}
\begin{array}{ccc}
J_1=\left(
\begin{array}{cc}
 a_{11} & 0  \\
 0 & a_{22}
\end{array}
\right), &
J_2=\left(
\begin{array}{cc}
 a_{11} & 1  \\
 -1 & a_{11}
\end{array}
\right), &
J_3=\left(
\begin{array}{cc}
 a_{11} & 1  \\
 0 & a_{11}
\end{array}
\right).
\end{array}
\end{equation}

The general solutions for $f$ and $g$ are listed as follows:
\begin{equation*}
  \begin{array}{ccc}
   \hline
    \text{Jordan form} & f & g \\
    \hline
    J_1 & \displaystyle{f_0 u^{-(4+\frac{4}{\gamma-1})}} & \displaystyle{f_1 u^{-\frac{4}{\gamma-1}}} \\
    J_2 & (f_0 y - f_1 z)\tau(y,z) & (f_0 z + f_1 y)\tau(y,z)\\
    J_3 & e^{-u}(f_0 u^{-4}+f_1 u^{-3}) & f_0 e^{-u}.\\
    \hline
  \end{array}
\end{equation*}
In this instance, $\tau(y,z) = e^{4 \alpha \arctan{\frac{z}{y}}}  (y^2 +z^2)^{-2}$, and $f_0\neq0,$ $f_1\neq0,$ $\gamma\neq1,$ $\alpha\neq1$ are constant.

Excluding reducible systems, the classes of functions $F$ and $G$ of equations (\ref{eq:detanl}) admitting a corresponding Lie group and the extension of the kernel of the admitted Lie algebras are obtained as seen in Table \ref{eq:xippn0}. The kernel of the admitted Lie algebras consists of the generator $X_1=\partial_x,$ which will be omitted on the list. The extension of the kernel is listed as follows:
\[
\begin{array}{lll}
  Y_2 =  2x  \partial_x  + y \partial_y  + z \partial_z, & Y_3 =  x (x  \partial_x  + y \partial_y  + z \partial_z), & \\
  Y_4 =  \gamma y \partial_y +  z \partial_z, &  Y_5 =  (\alpha y - z) \partial_y + (y + \alpha z) \partial_z,  & \\
  Y_6 =  (y + 4 z) \partial_y + z \partial_z, & Y_7 =  \cos{2x}  \partial_x  - \sin{2x}(y \partial_y  + z \partial_z),  &\\
  Y_8 =  \sin{2x}  \partial_x  + \cos{2x}(y \partial_y  + z \partial_z), & Y_9 =  e^{-2x} (\partial_x  - (y \partial_y  + z \partial_z )), & \\
  Y_{10} =  e^{2x} (\partial_x  + y \partial_y  + z \partial_z ). & &
\end{array}
\] The Lie algebras $Y_2,\, Y_3,\, Y_7,\, Y_8,\, Y_9,$ and $Y_{10}$ are associated with the coefficient $\xi$ and the Lie algebras $Y_4,\,Y_5,$ and $Y_6$ are related to the type of Jordan form of matrix $A$.

\subsection{Systems where all admitted generators have $\xi^{\prime\prime}=0$}
Note that the action of equivalence transformations coincides with the action of group automorphisms. For the direct approach, sometimes it is difficult to select out equivalent cases with respect to equivalence transformations. Fortunately, if the algebraic structure of the admitted Lie algebra is known, then using the algebraic approach aids in simplifying the group classification problem. Thus, for finding the group classification of systems of two autonomous nonlinear second-order ordinary differential equations with all admitted generators satisfying $\xi^{\prime\prime}=0,$ the one-dimensional optimal system of one parameter subgroups is utilized and is proceeded by the direct approach. Firstly, the two-step algorithm of Ovsiannikov \cite{bk:Ovsiannikov[1993opt]} is employed here, for which the optimal systems of subgroup and group invariant solutions are reconstructed. Invariant solutions are then substituted back into the determining equations (\ref{eq:detanl}) where the direct method is used to find all possible admitted Lie algebras.

Firstly returning to the analysis of the determining equations (\ref{eq:detanl}), since $\xi^{\prime\prime}=0$, it follows that $\xi = k_1 + k_2 x,$ where $k_1$ and $k_2$ are constant. This property of the coefficient forces $\zeta$ to become constant.\footnote{See in Appendix.}
The determining equations (\ref{eq:detanl}) are then reduced to
\begin{equation}\label{eq:detanl1}
\begin{array}{c}
3k_2 \mathbf{F} + (((A + k_2 E)\mathbf{y} + \mathbf{k})\cdot\nabla)\mathbf{F} - A \mathbf{F}=0

\end{array}
\end{equation}
with the following admitted generator
\begin{equation}\label{eq:gen}
X = 2 (k_1 + k_2x)\partial_x + (A\mathbf{y} + \mathbf{k})\cdot\nabla
\end{equation}
where the matrix $A$ is a vector $\mathbf{k},$  $k_1$ and $k_2$ are constant.
By rewriting (\ref{eq:gen}), the generator can be represented as
\begin{equation}\label{eq:gen1}
X=\sum\limits_{i=1}^8 c_i X_i
\end{equation}
where $c_i (i=1,\ldots8)$ are constant.
Corresponding to the constants $c_i (i=1,...,8)$, the basis operators of the  Lie algebra are as follows:
\begin{equation}
\begin{array}{cccc}
X_1=\partial_{x} & X_2=x \partial_{x} & X_3=\partial_{y} & X_4=\partial_{z}\\
X_5=y\partial_{y} & X_6=z \partial_{z} & X_7=z \partial_{y} & X_8=y \partial_{z}.\\
\end{array}
\end{equation}
From here, the one-dimensional optimal system of one parameter subgroups of the main group of system (\ref{eq:auto}) with $\xi^{\prime\prime}=0$ is constructed. The commutators of the basis operators are 
\begin{equation} 
\begin{array}{cccccccc}
\left[X_1,X_2\right] & = & X_1, & & & \left[X_5,X_7\right] & = & -X_7,\\
\left[X_3,X_5\right] & = & X_3, & & & \left[X_5,X_8\right] & = & X_8,\\
\left[X_3,X_8\right] & = & X_4, & & & \left[X_6,X_7\right] & = & X_7, \\
\left[X_4,X_6\right] & = & X_4, & & & \left[X_6,X_8\right] & = & -X_8, \\
\left[X_4,X_7\right] & = & X_3, & & & \left[X_7,X_8\right] & = & X_6-X_5.
\end{array}
\end{equation}
The following inner automorphisms $\textit{A}_i (i=1,...,8)$ of the above Lie algebra are found without difficulties:
\begin{equation}
\begin{array}{c}
A_1 : \hat{c}_1 = c_1 - a_1 c_2,\\
A_2 : \hat{c}_1 = e^{a_2} c_1, \\
A_3 : \hat{c}_3 = c_3 - a_3 c_5,\quad \hat{c}_4 = c_4 - a_3 c_8,\\
A_4 : \hat{c}_3 = c_3 - a_4 c_7,\quad \hat{c}_4 = c_4 - a_4 c_6,\\
A_5 : \hat{c}_3 = e^{a_5} c_3, \quad  \hat{c}_7 = e^{a_5} c_7 \quad \hat{c}_8 = e^{-a_5} c_8,\\
A_6 : \hat{c}_4 = e^{a_6} c_4, \quad  \hat{c}_7 = e^{-a_6} c_7 \quad \hat{c}_8 = e^{a_6} c_8,\\
A_7 : \hat{c}_3 = c_3 + a_7 c_4,\,\, \hat{c}_5 = c_5 + a_7 c_8, \,\,\hat{c}_6 = c_6 - a_7 c_8,\,\,
\\
\hspace*{20mm}\hat{c}_7 = c_7 - a_7^2 c_8 + a_7 c_6 - a_7 c_5,\\
A_8 : \hat{c}_4 = c_4 + a_8 c_3,\,\, \hat{c}_5 = c_5 - a_8 c_7, \,\, \hat{c}_6 = c_6 + a_8 c_7,\,\,
\\
\hspace*{20mm}\hat{c}_8 = c_8 - a_8^2 c_7 - a_8 c_6 + a_8 c_5.\\
\end{array}
\end{equation}
Note that $a_i \,(i=1,...,8)$ are the parameters on which the transformations of the group depend on. 
 Apart from these automorphisms, we have the following involutions:
\[
\begin{array}{c}
E_1:\, \bar{z}= - z |\quad \bar{c_4}= - c_4,\, \bar{c_7}= - c_7,\,\bar{c_8}= - c_8;\\
E_2:\, \bar{y}= - y |\quad \bar{c_3}= - c_3,\, \bar{c_7}= - c_7,\,\bar{c_8}= - c_8;\\
E_3:\, \bar{x}= - x |\quad \bar{c_1}= - c_1;\\
E_4:\, \bar{y}=  z,\, \bar{z}=  y |\quad \bar{c_3}= c_4,\, \bar{c_4}= c_3,\, \bar{c_5}= c_6,\, \bar{c_5}= c_5,\, \bar{c_7}= c_8,\,\bar{c_8}= c_7.\\
\end{array}
\]
We study the way in which the coefficients of equation (\ref{eq:gen1}) are changed under the action of inner automorphisms of the group above. Here and further on, only changeable coordinates of the generator are presented.
Looking closely at the commutators, the Lie algebra $L_8,$ which is composed of the generators $X_i$ $(i=1,...,8),$ can be split into 2 subalgebras $L_2 \oplus L_6 = $ ${\{X_1, X_2\}} \oplus{\{X_3, X_4, X_5, X_6, X_7, X_8\}}$.
Note also that $L_6$ can be decomposed further to $L_4\oplus I_2 = \{X_5, X_6, X_7, X_8\} \oplus \{X_3, X_4\}$, where $L_4$ makes up a 4-dimensional subalgebra and $I_2$ is ideal.

Let us first study the 4-dimensional subalgebra $L_4=\{X_5, X_6, X_7, X_8\}$. We consider this study here due to a misprint found in the classification of this Lie algebra in \cite{bk:PateraWinternitz[1977]}.  Now consider the operator $X$ of a one parameter subgroup of the form
\begin{equation}\label{eq:opsys}
X = c_5 X_5 +c_6 X_6+c_7 X_7+c_8 X_8.
\end{equation}
Automorphisms $A_5$ up to $A_8$ are made use of in order to find the one-dimensional optimal system of subalgebras of this Lie algebra.
From the automorphisms $A_5$ and $A_6,$ one can find the invariant $\bar{c_7}\bar{c_8}=c_7c_8,$ which leads one to consider the following cases:
\[
\begin{array}{cl}
(a)  & c_7c_8 > 0 \\
(b)  & c_7c_8 < 0\\
(c)  & c_7c_8 = 0.\\
\end{array}
\]
Then utilizing the invariant of $A_7$ and $A_8,$ which is $\bar{c_5}+\bar{c_6}=c_5 + c_6,$ one can obtain relations of $c_5$ and $c_6$. Upon further computations using automorphisms, it can be verified that for case $(a)$, the coefficients of equation (\ref{eq:opsys}) satisfy  $c_5 - c_6\neq 0,$ $c_7=0$ and $c_8=0.$  For case $(b)$,
it follows that $c_5 = c_6,$ $c_7=-1$ and $c_8=1.$ For case $(c)$ if $c_5\neq c_6$ then $c_7=0$ and $c_8=0.$ If $c_5 = c_6,$ then $c_7=1$ and $c_8=0.$ The involutions are also utilized. Hence, the following one-dimensional optimal system of subalgebras of the Lie algebra $L_4$ is obtained:
\begin{equation}\label{eq:opsys4}
\begin{array}{cl}
1. & X_5 + \alpha X_6 \,\,\text{ where } -1\leq \alpha\leq1\\
2. & \alpha (X_5 + X_6) + X_8 - X_7\,\, \text{ where } \alpha\geq0\\
3. & \beta (X_5 + X_6) + X_7 \text{ where } \beta = 0,1\\
4. & 0.\\
\end{array}
\end{equation}
Note that the $0$ element is considered on this list \cite{bk:Ovsiannikov[1993opt]}. There is a necessity to include this element on the list as when the direct sum $L_4  \oplus I_2$ is applied, more subalgebras of the Lie algebra $L_6$ may appear on the list.

\noindent \textbf{Remark}: As the action of the above automorphisms coincides with the action of the equivalence transformations, it is possible to get the optimal system of one-dimensional subalgebras of the Lie algebra $L_4$ using the latter. From the determining equations (\ref{eq:detanl1}) and the utilization of the equivalence transformation $\tilde{y}=Py,$ where $P$ is a nonsingular $2\times2$ matrix with constant entries, the matrix of coefficients of (\ref{eq:opsys}) \[
\left(
\begin{array}{cc}
 c_5 & c_7  \\
 c_8 & c_6
\end{array}
\right)
\] is reduced to one of its real-valued Jordan forms (\ref{eq:jor}). Looking closely at (\ref{eq:opsys4}),  subalgebra $1.$ coincides with Jordan matrix $J_1,$ subalgebra $2.$ coincides with Jordan matrix $J_2,$ and subalgebra $3.$ coincides with Jordan matrix $J_3.$

\subsubsection{Optimal system of subalgebras of the algebra $L_6 = \{X_3, X_4, X_5, X_6, X_7, X_8\}$}

After obtaining the one-dimensional optimal system (\ref{eq:opsys4}) of subalgebras of the Lie algebra $L_4 = \{X_5, X_6, X_7, X_8\}$, the next step is to combine $L_4$ with the ideal $I_2 = \{X_3, X_4\}.$ Here, again Ovsiannikov's two-step method \cite{bk:Ovsiannikov[1993opt]} is applied. Hence, for the study of the one-dimensional subalgebras of the Lie algebra $L_6,$ the study is reduced to analyzing the following elements:
\begin{equation}
\begin{array}{cl}
 1. & c_3 X_3 + c_4 X_4 + X_5 + \alpha X_6 \,\,\text{ where } -1\leq \alpha\leq1\\
 2. & c_3 X_3 + c_4 X_4 + \alpha (X_5 + X_6) + X_8 - X_7\,\, \text{ where } \alpha\geq0\\
 3. & c_3 X_3 + c_4 X_4 + \beta (X_5 + X_6) + X_7\,\, \text{ where } \beta = 0,1\\
 4. & c_3 X_3 + c_4 X_4.
\end{array}
\end{equation}

Using automorphisms $A_3$ and $A_4$ and the involutions, the list of one-dimensional subalgebras of the Lie algebra
$L_6 = \{X_3, X_4, X_5, X_6, X_7, X_8\}$ is obtained as follows:
\begin{equation}
\begin{array}{cl}
 1. & X_5 + \alpha X_6 \,\,\text{ where } -1\leq \alpha\leq1\\
 2. & X_4 + X_5\\
 3. & X_8 - X_7 \\
 4. & \beta X_3 + \alpha (X_5 + X_6) + X_8 - X_7\,\, \text{ where } \beta=-1,0,1,\,\alpha>0\\
 5. & \beta X_4 + X_7\,\, \text{ where } \beta = 0,1\\
 6. & X_5 + X_6 + X_7 \\
 7. & X_3\\
 8. & 0.
\end{array}
\end{equation}
Again, it is necessary to study the element $0$ of the subalgebras of the Lie algebra $L_6$ as this may generate additional elements when $L_6$ is combined with $L_2.$

\subsubsection{Optimal system of subalgebras of the algebra $L_8 = \{X_1, X_2, X_3, X_4, X_5, X_6, X_7, X_8\}$}
Combining $L_6$ with $L_2$ and keeping in mind that for autonomous systems $X_1$ is already admitted, the following elements comprise the list of one-dimensional subalgebras of the Lie algebra $L_8$:
\begin{equation}
\begin{array}{cl}
 1. & \gamma X_2 + X_5 + \alpha X_6 \,\,\text{ where } -1\leq \alpha\leq1\\
 2. & \gamma X_2 + X_4 + X_5\\
 3. & \gamma X_2 + X_8 - X_7 \\
 4. & \gamma X_2 + \beta X_3 + \alpha (X_5 + X_6) + X_8 - X_7\,\, \text{ where } \beta=-1,0,1,\,\alpha>0\\
 5. & \gamma X_2 + \beta X_4 + X_7 \,\,\text{ where } \beta = 0,1\\
 6. & \gamma X_2 + X_5 + X_6 + X_7 \\
 7. & \gamma X_2 + X_3\\
 8. &  X_2.
\end{array}
\end{equation}
Using this list of subalgebras, the next step is to obtain invariant solutions $F$ and $G$ of the determining equations (\ref{eq:detanl1}).
These functions are substituted into the determining equations (\ref{eq:detanl}), which are solved completely in order to find all other generators admitting equations (\ref{eq:auto}).

\subsubsection{Representations of systems of two nonlinear second-order ordinary differential equations with all generators having $\xi^{\prime\prime}=0$}
From (\ref{eq:gen1}), $c_i \,(i=1,\ldots,8)$ are the coefficients of the generator chosen from the above list of subalgebras. 
Only one subalgebra is presented in this paper as computations for the other subalgebras are done in a similar way.

\paragraph{Subalgebra 1. with the generator $\gamma X_2 + X_5 + \alpha X_6$  where $-1\leq \alpha\leq1$.}
For this case, the determining equations (\ref{eq:detanl1}) become
\[
\begin{array}{c}
   y F_y + \alpha z F_z - (2\gamma - 1)F = 0 \\
   y G_y + \alpha z G_z - (2\gamma - \alpha)G = 0.
\end{array}
\]
The general solution of these equations is
\[\begin{array}{cc}
  F(y,z)=f(u)y^{1-2\gamma} \text{ and }   & G(y,z)=g(u)y^{\alpha-2\gamma},
  \end{array}
\]
   where $u=y^\alpha/z.$ Notice that $f^{\prime}\neq0$ else it is equivalent to a reducible case. Substituting these functions to the determining equations (\ref{eq:detanl}), the following initial determining equations are obtained
\[
\begin{array}{l}
y^{2\alpha} a_{12} (\alpha u f^{\prime}  + (1-2\gamma)f - ug) \\
+ y^{\alpha +1}u ((\alpha a_{11} + (\alpha-1)\xi_1- a_{22})u f^{\prime} - 2(\gamma a_{11}+(\gamma-2)\xi_1)f )\\
+   y^{\alpha} \zeta_1 u (\alpha u f^{\prime} +  (1-2\gamma) f ) - y^2a_{21}u^3f^{\prime} -y \zeta_2 u^3f^{\prime} = 0,\\
y^{2\alpha}a_{12} ( \alpha u g^{\prime} + (\alpha-2\gamma)g)\\
+ y^{\alpha +1}u ( (\alpha a_{11} + (\alpha-1)\xi_1- a_{22})u g^{\prime}  - ((\alpha-2\gamma)a_{11} +(\alpha -2\gamma +3) \xi_1-a_{22} )g ) \\
+ y^{\alpha}u \zeta_1 (\alpha u g^{\prime} + (\alpha-2\gamma)g)- y^2a_{21}u (u^2g^{\prime} + f) -y \zeta_2 u^3 g^{\prime} = 0\\
\end{array}
\]
In order to split these determining equations, one needs to study relations between the powers of $y$. Thus, one needs to evaluate the following cases: (1) $\alpha = 0$, (2) $\alpha = \frac{1}{2}$, (3) $\alpha = 1$ and (4) $\alpha\neq0,\,\frac{1}{2},\,1$.
\begin{enumerate}
\item Consider when $\alpha =0.$ After splitting with respect to $y$, it can be verified that $a_{21}=0$ and one is left with the following determining equations
\begin{equation}\label{eq:c11}
\begin{array}{c}
( (1- 2 \gamma) (\zeta_1 u + a_{12}) )f - a_{12} u g=0,\\
2 ( \gamma a_{11} +(\gamma-2) \xi_1)f + (a_{22} + \xi_1 + \zeta_2 u) u f^{\prime}= 0, \\
\gamma (a_{12} + \zeta_1 u) g = 0,\\
( 2 \gamma a_{11} +(2 \gamma-3) \xi_1 + a_{22})g + (a_{22} + \xi_1 + \zeta_2u)u g^{\prime} =0.
\end{array}
\end{equation}
From the third equation, notice that if $\gamma=0$, $G$ becomes a function solely of $z$ and hence, this case is reducible. Thus, it follows that $a_{12}=0$ and  $\zeta_1=0.$ Dividing the second equation by $f^{\prime}$ and $u$, and differentiating it with respect to $u$ 2 times, one can study the following cases:
\begin{enumerate}
  \item $\left(\displaystyle{\frac{f}{u f^{\prime}}}\right)^{\prime\prime}\neq0$ and
  \item $\left(\displaystyle{\frac{f}{u f^{\prime}}}\right)^{\prime\prime}=0.$

  For the case when $\left(\displaystyle{\frac{f}{u f^{\prime}}}\right)^{\prime\prime}\neq0,$ it follows that $a_{11} =  \xi_1 \displaystyle{\frac{(2-\gamma)}{\gamma} }.$ Consequently, $\zeta_2=0$ and $a_{22} =- \xi_1.$  Substituting this into the determining equations (\ref{eq:c11}), no other extensions of the generator are obtained  apart from the studied subalgebra.

  For the case when $\left(\displaystyle{\frac{f}{u f^{\prime}}}\right)^{\prime\prime}=0,$ it follows that $\displaystyle{\frac{f}{uf^{\prime}}} = \kappa u + \beta.$ Furthermore, the general solution of this depends on $\beta.$ Thus, one needs to study whether $\beta\neq0$ or $\beta=0.$
    \begin{enumerate}
      \item For the case when $\beta\neq0,$ the general equation for $f$ (with a possible shift) is $f_0 \left(\displaystyle{\frac{1}{u}} \right)^{\beta}.$ Substituting this into the determining equations, one gets $a_{22} = \displaystyle{\frac{(2\gamma -\beta-4)\xi_1 +2\gamma a_{11}}{\beta}}$ and $\zeta_2 = \displaystyle{\frac{2\kappa(\gamma a_{11}+ (\gamma-2) \xi_1)}{\beta}}.$ Consequently, the general solution for $g$ is $g_0 \left(\displaystyle{\frac{1}{u}} \right)^{\beta+1}.$ The extension $\beta X_5 + 2\gamma X_6$ is obtained along with the studied subalgebra.
      \item For the case when $\beta=0,$ it follows that $\kappa\neq0.$ Hence, the general equation for $f$ is $f_0 e^{\kappa/u}.$ Substituting this into the determining equations, one obtains $\zeta_2 = \displaystyle{\frac{2(\gamma a_{11}+ (\gamma-2) \xi_1)}{\kappa}}$ and $a_{22} = -\xi_1.$ Consequently, the general solution for $g$ is $g_0 e^{\kappa/u}.$ The extension $\kappa X_5 + 2 \gamma X_4$ is obtained aside the studied subalgebra. 
    \end{enumerate}
\end{enumerate}

\item For the case when $\alpha=\frac{1}{2},$ after splitting with respect to $y,$ it follows that $a_{21}=0.$ Also, since $(1- 4\gamma)g + ug^{\prime}=0$ leads to a reducible case it then follows that $\zeta_1=0.$ The remaining determining equations are
    \[
    \begin{array}{c}
    2 a_{12} (1 - 2 \gamma) f - 2 a_{12} u g + (a_{12} - 2 \zeta_2 u^2) u f^{\prime} = 0,\\
    4 ((2- \gamma) \xi_1 - \gamma a_{11})f  + (a_{11} - 2 a_{22} - \xi_1) u f^{\prime} =0,\\
    (1 - 4 \gamma) a_{12} g  + (a_{12} - 2 \zeta_2 u^2) u g^{\prime}=0,\\
    ( (1- 4 \gamma) a_{11} +(7 - 4 \gamma) \xi_1 - 2 a_{22} ) g + (a_{11} - 2 a_{22} - \xi_1) u g^{\prime}=0.
    \end{array}
   \]
   Dividing the fourth equation by $g$ (as it is nonzero) and differentiating it with respect to $u$, one is left to study the following cases:
   \begin{enumerate}
     \item $\left(\displaystyle{\frac{u g^{\prime}}{g}}\right)^{\prime}\neq0$ and
     \item $\left(\displaystyle{\frac{u g^{\prime}}{g}}\right)^{\prime}=0.$

     Consider $\left(\displaystyle{\frac{u g^{\prime}}{g}}\right)^{\prime}\neq0.$ It follows that $a_{11} = 2 a_{22} + \xi_1.$ If $\gamma=0$ then $\xi_1=0$, but if $\gamma\neq0$ then $a_{22}=\xi_1 \left(\displaystyle{\frac{1-\gamma}{\gamma}}\right).$
     From the third equation, one needs to study the following cases:
     \begin{enumerate}
       \item the case where there exists a generator with $a_{12}\neq0,$ and
       \item the case for which all generators have $a_{12}=0.$

       If there exists a generator with $a_{12}\neq0,$ then $g$ satisfies the form $(1-4\gamma)g + (1-\beta u^2)u g^{\prime}=0.$ Notice that $\beta=0$ is reducible. Hence, $\beta\neq0.$ Without loss of generality, one can assume that $\beta=1.$ Then the general solution of $g$ is $g_0 \left(1 -\displaystyle{\frac{1}{u^2}}\right)^{\tilde{\gamma}}$ where $\tilde{\gamma}=\displaystyle{\frac{1-4\gamma}{2}}\neq0$ (if $\tilde{\gamma}=0$, the case is reducible). Substituting this into the determining equations, we obtain $\zeta_2=\displaystyle{\frac{a_{12}}{2}}.$ It follows that $f=\phi(u) \left(1 -\displaystyle{\frac{1}{u^2}}\right)^{\tilde{\gamma}+(1/2)},$ where $\phi = f_0 - 2 g_0 \left( \displaystyle{\frac{1}{(u^2-1)^{(1/2)}}}\right).$  Here, the extension $X_4 + 2 X_7$ is obtained besides the studied subalgebra.

       For the case where all generators have $a_{12}=0,$ it follows that $\zeta_2=0.$ No other extensions are obtained for this case.
     \end{enumerate}

     Consider $\left(\displaystyle{\frac{u g^{\prime}}{g}}\right)^{\prime}=0.$ The general solution for this is $g=g_0 u^{\kappa}$. Substituting this into the determining equations, one obtains that $a_{12}=0$ and $\zeta_2=0.$ Consequently, the form of $f$ either satisfies $(\kappa+1) f -u f^{\prime} =0$ or not. If it is satisfied, then the general solution is $f=f_0 u^{\kappa+1}.$ Moreover, $a_{22} = (\kappa - 4\gamma + 1)(a_{11} - \xi_1) + 8\xi_1.$ Here, the extension $(\kappa+1)X_2 + 2 X_6$ is obtained apart from the studied subalgebra. If $f$ does not satisfy $(\kappa+1) f - f^{\prime} u =0,$ then no extensions are obtained other than the studied subalgebra.

   \end{enumerate}

\item For the case when $\alpha=1,$ the determining equations after splitting with respect to $y$ are as follows
\[
\begin{array}{c}
  (1 - 2 \gamma)\zeta_1 f + ( \zeta_1 - \zeta_2 u)u f^{\prime} =0, \\
  ((1- 2 \gamma) a_{12} +((4- 2 \gamma) \xi_1 u  - 2 \gamma a_{11} u))f-  a_{12} u g \\
  + ((a_{11}- a_{22}) u + a_{12} - a_{21} u^2 ) u f^{\prime}=0,\\
  (1 - 2 \gamma) \zeta_1 g + ( \zeta_1 - \zeta_2 u)u g^{\prime}  =0,\\
   - a_{21} u f + g ( (1- 2 \gamma) a_{11} u +(1- 2 \gamma) a_{12} +(4- 2 \gamma) \xi_1 u  - a_{22} u ) \\
  + g^{\prime} u ((a_{11}- a_{22}) u + a_{12} - a_{21} u^2 ) =0.\\
\end{array}
\]
From the first and third equations, one can study the following 2 cases:
\begin{enumerate}
  \item $fg^{\prime}-gf^{\prime}=0,$ and
  \item $fg^{\prime}-gf^{\prime}\neq0.$

  Notice that when $fg^{\prime}-gf^{\prime}=0,$ then $g=g_0 f$. This is a reducible case. 
  Hence, we consider only when $fg^{\prime}-gf^{\prime}\neq0.$ It follows that $\zeta_1=\zeta_2=0.$
  From here, one can assume that $g=\phi(u) f$ (as $f$ is nonzero), where $\phi^{\prime}\neq0$.  If it is assumed further that $\phi = \psi(u) + 1/u,$ then the determining equations are reduced as follows:
  \begin{equation}\label{eq:c13}
    \begin{array}{c}
       (2 ( - \gamma a_{11} u + (2- \gamma) \xi_1 u) - (\gamma  + \psi u) a_{12})f \\
      +  ((a_{11}-a_{22}) u + a_{12} - a_{21} u^2 ) u f^{\prime}=0,\\
      (a_{11}-a_{22}) u + a_{12} - a_{21} u^2 ) \psi^{\prime} +  a_{12}\psi^2 + (a_{11} + 2 a_{12} u^{-1}   - a_{22})\psi = 0.
    \end{array}
  \end{equation}
  These equations lead one to study the two cases where:
  \begin{enumerate}
    \item there exists at least one generator with $a_{12} \neq 0,$ and
    \item where all generators have $a_{12} = 0.$

    For the case where there exists at least one generator with $a_{12} \neq 0,$ it follows that $\psi(u) = -\displaystyle{\frac{\kappa u^2+\lambda u +\beta}{u (\beta - \psi_0 u)}},$ where $\beta\neq0,\,\psi_0\neq0,\lambda,\,\kappa$ are constant. Without loss of generality, it is assumed further that $\beta=1$. Consequently, we obtain $a_{11} = \lambda a_{12}+ a_{22}$ and $a_{21} = -\kappa a_{12}$. Substituting this into the remaining determining equations gives the solution for $f$ which depends on the following three cases:
    \begin{enumerate}
      \item $4 \kappa - \lambda^2 > 0,$
      \item $4 \kappa - \lambda^2 < 0,$ and
      \item $4 \kappa - \lambda^2 = 0.$

      For the case where $4 \kappa - \lambda^2 > 0,$ it is assumed that $4 \kappa - \lambda^2 = p^2,\,p\neq0.$ The solution for $f$ is
      \[f_0 \displaystyle{\frac{(1-\psi_0 u) u^{2\gamma-1}}{(\kappa u^2+\lambda u +1)^{\gamma}}} e^{\left(\displaystyle{\frac{(2\lambda\gamma-4\mu)}{p} \arctan\left(\frac{\lambda+2\kappa u}{p}\right)}\right)}\]
      where $\mu$ is constant.

      For the case where $4 \kappa - \lambda^2 < 0,$ it is assumed that $4 \kappa - \lambda^2 = -p^2,\,p\neq0.$ The solution for $f$ is
      \[f_0 \displaystyle{\frac{(1-\psi_0 u) u^{2\gamma-1}}{(\kappa u^2+\lambda u +1)^{\gamma}}} \displaystyle{\left(\frac{2\kappa u +\lambda-p}{2\beta\kappa u +\lambda+p}\right)^{\frac{\lambda\gamma-2\mu}{p}}} \]
      where $\mu_0$ is constant.

      For the case where $4 \kappa - \lambda^2 = 0,$ it follows that
      \[f = f_0 \displaystyle{\frac{(1-\psi_0 u) u^{2\gamma-1}}{(\kappa u^2+\lambda u +1)^{\gamma}}} e^{\left(-\displaystyle{\frac{4(\gamma + \mu u)}{\lambda u + 2} }\right)}\] where $\mu$ is constant.

       If  $\gamma\neq0,$ then $a_{22}= \displaystyle{\frac{(2-\gamma) \xi_1 -\mu a_{12}}{\gamma}}$. If $\gamma=0,$ then $\xi_1= \displaystyle{\frac{\mu a_{12}}{2}}$. Here, the extension $(\lambda\gamma - \mu) X_5 -  \mu X_6 + \gamma X_7 - \kappa\gamma X_8$ is obtained apart from the studied subalgebra.


    \end{enumerate}

    For the case where all generators have $a_{12} = 0,$ the determining equations are reduced to
    \begin{equation}\label{eq:c112}
    \begin{array}{c}
      2( (2-\gamma)\xi_1 - \gamma a_{11}) f + (a_{11} - a_{22} -a_{21}u) u f^{\prime} = 0, \\
      ((4-2\gamma) \xi_1 -a_{22} +(1-2\gamma) a_{11}) g - a_{21}f  \\
      + (a_{11} - a_{22} -a_{21}u) u g^{\prime} = 0.
    \end{array}
    \end{equation}
    Dividing the first equation (\ref{eq:c112}) by $uf^{\prime}$ and differentiating this equation with respect to $u$ twice, leads to the study of the following sub-cases:
    \begin{enumerate}
      \item $\left(\displaystyle{\frac{f}{uf^{\prime}}}\right)^{\prime\prime}\neq0,$ and
      \item $\left(\displaystyle{\frac{f}{uf^{\prime}}}\right)^{\prime\prime}=0.$

      If $\left(\displaystyle{\frac{f}{uf^{\prime}}}\right)^{\prime\prime}\neq0,$ then it follows that if $\gamma\neq0$ then $a_{11} =  \xi_1 \displaystyle{\frac{2-\gamma}{\gamma}}$, $a_{22} = \xi_1 \displaystyle{\frac{2-\gamma}{\gamma}}$ and $a_{21}=0.$ If $\gamma=0$ then $\xi_1=0,$ $a_{22} = a_{11}$ and $a_{21}=0.$ For both cases, no extensions are obtained apart from the studied subalgebra. 

      If $\left(\displaystyle{\frac{f}{f^{\prime}u}}\right)^{\prime\prime}=0,$ then the general solution for $f$ is $f_0 \left(\displaystyle{\frac{u}{1+u}}\right)^{\kappa},$ where $\kappa\neq0$ (else it is reducible). Substituting this into the determining equations (\ref{eq:c112}), one obtains that $a_{21}= 2\left(\displaystyle{\frac{-\gamma a_{11} +(2-\gamma)\xi_1}{\kappa}}\right) $ and $a_{22}= \displaystyle{\frac{(\kappa-2\gamma) a_{11} +(4-2\gamma)\xi_1}{\kappa}}.$ Substituting this into the remaining determining equation, one finds that $g$  satisfies $g^{\prime} u (1+u) +(1-\kappa) g +f =0.$ The general solution of this is $g=\left(g_0 - f_0 \displaystyle{\frac{u}{u+1}}\right) \left(\displaystyle{\frac{u}{u+1}}\right)^{\kappa-1}.$ The extension $\kappa X_2 + 2(X_6+X_8)$ is obtained aside from the studied subalgebra. 
    \end{enumerate}

  \end{enumerate}

\end{enumerate}

\item For the case where $\alpha\neq0,\,\frac{1}{2},\,1,$ the determining equations are split with respect to $y$. Since $f^{\prime}\neq0$, it follows that $\zeta_2=0$ and $a_{21}=0.$ Notice also that since $\alpha u g^{\prime}  + (\alpha -2\gamma)g =0$ leads to a degenerate case, then $\zeta_1=0$ and $a_{12}=0.$ The remaining determining equations become
\begin{equation}\label{eq:c14}
\begin{array}{c}
(\alpha a_{11} + (\alpha - 1) \xi_1 - a_{22} )u f^{\prime} +  (- 2\gamma a_{11} + (4- 2\gamma) \xi_1)f =0, \\
(\alpha a_{11} + (\alpha - 1) \xi_1 - a_{22} )u g^{\prime} +  ((\alpha- 2\gamma) a_{11} + (\alpha - 2\gamma + 3) \xi_1  - a_{22}) g =0.
\end{array}
\end{equation}
Dividing the first equation by $f$ (as it is nonzero) and differentiating with respect to $u$, it can be observed that there is need to study the following cases:
\begin{enumerate}
  \item $\left(\displaystyle{\frac{uf^{\prime}}{f}}\right)^{\prime}\neq0$ and
  \item $\left(\displaystyle{\frac{uf^{\prime}}{f}}\right)^{\prime}=0.$

  For the case with $\left(\displaystyle{\frac{uf^{\prime}}{f}}\right)^{\prime}\neq0$, it follows that $a_{22} =\alpha a_{11} + (\alpha -1)\xi_1.$  Substituting this into the remaining determining equations (\ref{eq:c14}), we find that if $\gamma\neq0$ then $a_{11}=\displaystyle{\frac{2-\gamma}{\gamma}}\xi_1$, and if $\gamma=0$ then $\xi_1=0.$ Substituting all these, no other extensions of the generator is found other than the studied subalgebra. 

  For the case where $\left(\displaystyle{\frac{uf^{\prime}}{f}}\right)^{\prime}=0$, the general solution for $f$ is $f_0 u^{\kappa}$, where $\kappa\neq0.$ Substituting this function into the determining equations, it follows that $a_{22} = \displaystyle{\frac{ (\kappa \alpha -2 \gamma)a_{11} + (\kappa \alpha -\kappa -2\gamma +4)\xi_1}{\kappa}}.$ This leads us to study the two cases, that is, if $g^{\prime}u  - g (\kappa - 1)=0$ or $g^{\prime}u  - g (\kappa - 1)\neq0.$ If $g^{\prime}u  - g (\kappa - 1)=0$, the general solution for $g$ is $g_0 u^{\kappa-1}.$  Another extension of the generator apart from the studied subalgebra is found, that is $\kappa X_2 + 2 X_6$. For the case where $g^{\prime}u  - g (\kappa - 1)\neq0$, if $\gamma\neq0$ then $a_{11}=\displaystyle{\frac{2-\gamma}{\gamma}}\xi_1$, and if $\gamma=0$ then $\xi_1=0.$ Again, these lead to a generator with only the studied subalgebra as its extension. 

\end{enumerate}

\end{enumerate}

The complete representative classes for the autonomous system with all admitted generators having $\xi^{\prime\prime}=0$ is listed in Tables \ref{eq:xipp0} and \ref{eq:xipp0d}.

\section{Conclusion}
A complete group classification of the systems of two autonomous nonlinear second-order ordinary differential equations of the form $\bf{y^{\prime\prime}}=\bf{F(y)}$ excluding the systems which are equivalent to linear systems and the degenerate case were presented using both the direct and algebraic approach. The important thing in this study is that the analysis of the determining equations where split into two cases: 1) the case where at least one admitted generator has $\xi^{\prime\prime}\neq0$ and 2) the case where all admitted generators have $\xi^{\prime\prime}=0.$ The first was analyzed through the direct approach while the latter was analyzed using one-dimensional optimal system of subalgebras followed by the direct approach. For the direct approach, all possible Lie algebras were found with the aid of the equivalence transformations applied on the determining equations. As for the algebraic approach, the study was reduced to the analysis of relations between constants of the generator with its corresponding basis operators. The obtained classification is summarized on Tables \ref{eq:xippn0}, \ref{eq:xipp0} and \ref{eq:xipp0d}. It is highly likely that the same methods shown in this paper are applicable to the group classification of systems of two nonlinear second-order ordinary differential equations, which will be next goal for further studies. In addition, it is also believed that this can be extended to systems in more general cases.

\section*{Acknowledgements}
The authors thank an anonymous referee for his valuable remarks to improve this paper.
GFO also thanks SUT-PhD Scholarship for ASEAN for the financial support during her study in Suranaree University of Technology.



\begin{thebibliography}{10}
\expandafter\ifx\csname url\endcsname\relax
  \def\url#1{\texttt{#1}}\fi
\expandafter\ifx\csname urlprefix\endcsname\relax\def\urlprefix{URL }\fi
\expandafter\ifx\csname href\endcsname\relax
  \def\href#1#2{#2} \def\path#1{#1}\fi

\bibitem{bk:Lie[1883]}
S.~Lie, {K}lassifikation und {I}ntegration von gew\"ohnlichen
  {D}ifferentialgleichungen zwischen $x,y$, die eine {G}ruppe von
  {T}ransformationen gestatten. {III}, Archiv for Matematik og Naturvidenskab
  8~(4) (1883) 371--427, reprinted in Lie's Gessammelte Abhandlungen, 1924, 5,
  paper XIY, pp. 362--427.

\bibitem{bk:Lie[1891b]}
S.~Lie, Vorlesungen \"{u}ber {D}ifferentialgleichungen mit bekannten
  infinitesimalen {T}ransformationen, B.G.Teubner, Leipzig, 1891, bearbeitet
  und herausgegeben von Dr. G.Scheffers.

\bibitem{bk:Lie[1884]}
S.~Lie, {K}lassifikation und {I}ntegration von gew\"ohnlichen
  {D}ifferentialgleichungen zwischen $x,y$, die eine {G}ruppe von
  {T}ransformationen gestatten. {IV}, Archiv for Matematik og Naturvidenskab
  8~(4) (1884) 431--448, reprinted in Lie's Gessammelte Abhandlungen, 1924, 5,
  paper XYI, pp. 432--446.

\bibitem{bk:Lie[1881]}
S.~Lie, \"{U}ber die {I}ntegration durch bestimmte {I}ntegrale von einer
  {K}lasse linearer partieller {D}ifferentialgleichungen, Arch. for Math. 6
  (1881) 328.

\bibitem{bk:ovsiannikov[1958]}
L.~V. Ovsiannikov, Groups and group-invariant solutions of partial differential
  equations, Dokl. AS USSR 118~(3) (1958) 439--442.

\bibitem{bk:Ovsiannikov[1978]}
L.~V. Ovsiannikov, Group analysis of differential equations, Nauka, Moscow,
  1978, {E}nglish translation, {A}mes, {W}.{F}., Ed., published by Academic
  Press, New York, 1982.

\bibitem{bk:MahomedLeach[1989]}
F.~M. Mahomed, P.~G.~L. Leach, {L}ie algebras associated with scalar
  second-order ordinary differential equations, Journal of Mathematical Physics
  30~(12) (1989) 2770--2777.

\bibitem{bk:GonzalezKamranOlver[1992a]}
A.~Gonzalez-Lopez, N.~Kamran, P.~J. Olver, Lie algebras of differential
  operators in two complex variables, American Journal of Mathematics 114
  (1992) 1163--1185.

\bibitem{bk:PopovychKunzingerEshraghi[2010]}
R.~O. Popovych, M.~Kunzinger, H.~Eshraghi, Admissible transformations and
  normalized classes of nonlinear {S}chr\"odinger equations, Acta Appl. Math.
  109 (2010) 315--359.

\bibitem{bk:GrigorievMeleshkoSuriyawichitseranee[2013]}
Y.~N.Grigoriev, S.~V. Meleshko, A.~Suriyawichitseranee, On the equation for the
  power moment generating function of the {B}oltzmann equation. group
  classification with respect to a source function, in: O.~Vaneeva,
  C.~Sophocleous, R.~Popovych, P.~Leach, V.~Boyko, P.~Damianou (Eds.), Group
  Analysis of Differential Equations \& Integrable Systems, University of
  Cyprus, Nicosia, 2013, pp. 98--110.

\bibitem{bk:Ovsiannikov[2004]}
L.~V. Ovsiannikov, Group classification of equation of the form
  $y^{\prime\prime}=f(x,y)$, Journal of Applied Mechanics and Technical Physics
  45~(2) (2004) 153--157.

\bibitem{bk:Phauk[2013]}
S.~Phauk, Group classification of second-order ordinary differential equations
  in the form of a cubic polynomial in the first-order derivative, mS Thesis at
  School of Mathematics, Suranaree University of Technology, Thailand (2013).

\bibitem{bk:MahomedLeach[1990]}
F.~M. Mahomed, P.~G.~L. Leach, Symmetry {L}ie algebras of $n$th order ordinary
  differential equations, Journal of Mathematical Analysis and Applications
  151~(12) (1990) 80--107.

\bibitem{bk:HandbookLie_v3}
N.~H. Ibragimov (Ed.), {CRC} Handbook of {L}ie Group Analysis of Differential
  Equations, Vol.~3, CRC Press, Boca Raton, 1996.

\bibitem{bk:Gorringeleach[1988]}
V.~Gorringe, P.~Leach, Lie point symmetries for systems of $2$nd order linear
  ordinary differential equations, Quaestiones Mathematicae 1 (1988) 95--117.

\bibitem{bk:WafoMahomed[2000]}
C.~{Wafo Soh}, F.~M. Mahomed, Symmetry breaking for a system of two linear
  second-order ordinary differential equations, Nonlinear Dynamics 22 (2000)
  121--133.

\bibitem{bk:BoyPopoSha[2012]}
V.~M. Boyko, R.~O. Popovych, N.~M. Shapoval, Lie symmetries of systems of
  second-order linear ordinary differential equations with constant
  coefficients, Journal of Mathematical Analysis and Applications 397 (2012)
  434--440.

\bibitem{bk:WafoSoh[2010]}
C.~{Wafo Soh}, Symmetry breaking of systems of linear second-order differential
  equations with constant coefficients, Communications in Nonlinear Science and
  Numerical Simulations 15 (2010) 139--143.

\bibitem{bk:Meleshko[2011]}
S.~V. Meleshko, Comment on "{S}ymetry breaking of systems of linear
  second-order ordinary differential equations with constant coefficients",
  Communications in Nonlinear Science and Numerical Simulations 16 (2011)
  3447--3450.

\bibitem{bk:Campoamor-Stursberg[2011]}
R.~Campoamor-Stursberg, Systems of second-order linear ode's with constant
  coefficients and their symmetries, Communications in Nonlinear Science and
  Numerical Simulations 16 (2011) 3015--3023.

\bibitem{bk:Campoamor-Stursberg[2012]}
R.~Campoamor-Stursberg, Systems of second-order linear ode's with constant
  coefficients and their symmetries. {II}, Communications in Nonlinear Science
  and Numerical Simulations 17 (2012) 1178--1193.

\bibitem{bk:MoyoMeleshkoOguis[2013]}
S.~Moyo, S.~V. Meleshko, G.~F. Oguis, Complete group classification of systems
  of two linear second-order ordinary differential equations, Communications in
  Nonlinear Science and Numerical Simulation 18~(11) (2013) 2972--2983.

\bibitem{bk:MkhizeMoyoMeleshko[2014]}
T.~G. Mkhize, S.~Moyo, S.~V. Meleshko, Complete group classification of systems
  of two linear second-order ordinary differential equations. {A}lgebraic
  approach, Mathematical Methods in the Applied Sciences 38 (2015) 1824--1837.

\bibitem{bk:MoyoMeleshkoOguis[2014]}
S.~V. Meleshko, S.~Moyo, G.~F. Oguis, On the group classification of systems of
  two linear second-order with constant coefficients, Journal of Mathematical
  Analysis and Applications 410 (2014) 341--347.

\bibitem{bk:MeleshkoMoyo[2014]}
S.~V. Meleshko, S.~Moyo, On the study of the general group classification of
  systems of linear second-order ordinary differential equations,
  Communications in Nonlinear Science and Numerical Simulation 22 (2015)
  1002--1016.

\bibitem{bk:Ovsiannikov[1993opt]}
L.~V. Ovsiannikov, On optimal system of subalgebras, Docl. RAS 333~(6) (1993)
  702--704.

\bibitem{bk:PateraWinternitz[1977]}
J.~Patera, P.~Winternitz, Subalgebras of real three- and four-dimensional lie
  algebras, Journal of Math. Phys. 18~(7) (1977) 1449--1455.

\end{thebibliography}

\newpage
\begin{center}
\begin{table}[!htb]
\protect\protect\caption{Group classification of systems admitting at least one generator with $\xi^{\prime\prime}\neq0$.}
\label{eq:xippn0}
\begin{tabular}{ccccc}
\hline
$F$ & $G$ & $\kappa$ & Extension of Kernel\\
\hline
$\kappa y +\frac{f_0 y}{z^{4}} ({{z}\over{y}})^{-4/(\gamma-1)}$ & $\kappa z + \frac{f_1}{z^3} ({{z}\over{y}})^{-4/(\gamma-1)}$ & $0$  & $Y_2, \,Y_3,\,Y_4$ \\
 &  & $-1$ & $Y_7,\, Y_8,\,Y_4$ \\
 &  & $1$  &  $Y_9,\, Y_{10},\,Y_4$ \\
$\kappa y + (f_0 y - f_1 z)\tau(y,z)$ & $\kappa z + (f_0 z + f_1 y)\tau(y,z)$ & $0$  &  $Y_2,\, Y_3,\,Y_5$ \\
&  & $-1$   &  $Y_7,\, Y_8,\,Y_5$ \\
& & $1$  &   $Y_9,\, Y_{10},\,Y_5$ \\
$\kappa y + {e^{{y}\over{z}} z^{-4} (f_0 y + f_1 z)}$ & $\kappa z + f_0 z^{-3} e^{{y}\over{z}}$ &  $0$ & $Y_2,\, Y_3,\,Y_6$ \\
&  & $-1$   & $Y_7,\, Y_8,\,Y_6$ \\
&  & $1$   & $Y_9,\, Y_{10},\,Y_6$ \\
\hline
\end{tabular}
\end{table}
\end{center}

\begin{landscape}
\begin{center}
\begin{table}[!htb]
\protect\protect\caption{Group classification of systems admitting all generator with $\xi^{\prime\prime}=0$. Here we have $\theta_1(u,v)=(\cos(u)f(v)+\sin(u)g(v))$,\newline $\theta_2(u,v)=\sin(u)f(v)-\cos(u)g(v)$, $\chi_1(\alpha)=\displaystyle{\frac{\alpha}{\alpha^2+1}}$
and $\chi_2(\alpha)=\displaystyle{\frac{1}{\alpha^2+1}}$.}
\label{eq:xipp0}
\begin{tabular}{cccc}
\hline
 $F$ & $G$ & Relations & Extension of Kernel\\
 \hline
 $f(u)y^{(1-2\gamma)}$ & $g(u)y^{(\alpha-2\gamma)}$ & $u=\frac{y^\alpha}{z}$ \,$-1\leq\alpha\leq1$ &  $\gamma X_2 + X_5 + \alpha X_6$ \\
 $f(u)y^{(1-2\gamma)}$ & $g(u)y^{(-2\gamma)}$ & $u=ye^{-z}$ &  $\gamma X_2 + X_4 + X_5$ \\
 $e^{-2\gamma u}\theta_1(u,v)$ & $-e^{-2\gamma u}\theta_2(u,v)$ & $y=v\cos(u),\,\,z=v\sin(u)$ & $\gamma X_2 - X_7 + X_8$\\
 $e^{(\alpha-2\gamma) u}\theta_1(u,v)$ & $e^{(\alpha-2\gamma) u}\theta_2(u,v)$ & $y=ve^{\alpha u}\cos(u)+\chi_1(\alpha)$ & \\
  &   &  $z=ve^{\alpha u}\sin(u)-\chi_2(\alpha),\,\alpha>0$ & $\gamma X_2 - X_3  + \alpha(X_5+X_6) -  X_7 + X_8$\\
 $e^{(\alpha-2\gamma) u}\theta_1(u,v)$ & $e^{(\alpha-2\gamma) u}\theta_2(u,v)$ & $y=ve^{\alpha u}\cos(u)-\chi_1(\alpha)$ & \\
  &   & $z=ve^{\alpha u}\sin(u)+\chi_2(\alpha),\,\alpha>0$ & $\gamma X_2 + X_3  + \alpha(X_5+X_6) -  X_7 + X_8$\\
 $e^{(\alpha-2\gamma) u}\theta_1(u,v)$ & $e^{(\alpha-2\gamma) u}\theta_2(u,v)$ & $y=ve^{\alpha u}\cos(u)$ & \\
  &   &  $z=ve^{\alpha u}\sin(u),\,\alpha>0$ & $\gamma X_2 + \alpha(X_5+X_6) -  X_7 + X_8$\\
 $(g(v)u+f(v))e^{(-2\gamma u)}$ & $g(v) e^{(-2\gamma u)}$ & $y=uv,\,\, z=v$ & $\gamma X_2+X_7$\\
 $(g(u)z+f(u))e^{(-2\gamma z)}$ & $g(u) e^{(-2\gamma z)}$ & $u=z^2-2y$ & $\gamma X_2+X_4+X_7$\\
 $((y/z)g(u)+f(u))e^{((1-2\gamma)(y/z))}$ & $g(u)e^{((1-2\gamma)(y/z))}$ & $u=z e^{-y/z}$ & $\gamma X_2 + X_5 + X_6 + X_7 $\\
 $f(z)e^{-2\gamma y}$ & $g(z)e^{-2\gamma y}$ & & $\gamma X_2 +X_3$\\
 \hline
\end{tabular}
\end{table}
\end{center}
\pagebreak

\setlength{\LTcapwidth}{\linewidth}

{\scriptsize
\begin{center}
\begin{longtable}{cccc}

 \caption{Group classification of systems admitting all generator with $\xi^{\prime\prime}=0$. Here, $f_0,\,g_0,\,\phi_0,\,\phi_1,\,\alpha,\,\beta,\,\kappa,\,\mu_0,\,\lambda$ and $\gamma$ are constant. 
}\\
\hline
 \endfirsthead

 \multicolumn{4}{c}%
{{\bfseries \tablename\ \thetable{} -- continued from previous page}} \\
 \hline
\endhead

\hline
\multicolumn{4}{r}{{\tablename\ \thetable{} -- continued on next page}} \\
\hline
\endfoot

\endlastfoot

 \multicolumn{4}{c}%
 {{\bfseries Subalgebra 1. $\gamma X_2 + X_5$ }} \\

 $F$ & $G$ & Relations & Additional Extension of Kernel \\
 $f_0 z^{\beta}y^{1+\tilde{\gamma}}$ & $g_0z^{\beta+1}y^{\tilde{\gamma}}$ & $\tilde{\gamma}=-2\gamma\neq0,\, \beta\neq0$  &   $\beta X_5 - \tilde{\gamma} X_6$ \\

 $f_0 y^{1+\tilde{\gamma}}e^{\kappa z}$ & $g_0 y^{\tilde{\gamma}}e^{\kappa z}$ &$\tilde{\gamma}=-2\gamma\neq0,\,\kappa\neq0$ & $ \kappa X_5 -\tilde{ \gamma} X_4$ \\

 \hline\\
 \hline\\
 \multicolumn{4}{c}%
 {{\bfseries Subalgebra 1. $\gamma X_2 + X_5 +\displaystyle{\frac{1}{2}}X_6$ }} \\


 $F$ & $G$ & Relations & Additional Extension of Kernel \\
 $\phi(u) (y-z^2)^{\tilde{\gamma}}$ & $g_0 (y-z^2)^{\tilde{\gamma}}$  & $\tilde{\gamma}=\displaystyle{\frac{1-4\gamma}{2}}\neq0$  & \\
    &    & $\phi = f_0 (y-z^2)^{1/2} + 2 g_0z$   & $  X_4 + 2 X_7$ \\
 $f_0 z^{-(\kappa+1)}y^{\tilde{\gamma}+1}$ & $g_0 z^{-\kappa}y^{\tilde{\gamma}}$ 
    &    $\kappa+1\neq0,\,\tilde{\gamma}=\displaystyle{\frac{k+1-4\gamma}{2}}\neq0$ & $ (\kappa+1) X_2 + 2 X_6$\\
 \hline\\
 \hline\\
 \multicolumn{4}{c}%
 {{\bfseries Subalgebra 1. $\gamma X_2 + X_5 + X_6$ }} \\
 $F$ & $G$ & Relations & Additional Extension of Kernel \\

 $ f_0 \displaystyle{\frac{z-\alpha y}{(z^2+\lambda yz +\kappa y^2)^{\gamma}}} \psi_i(y,z) $ & $- f_0 \displaystyle{\frac{(\kappa y + (\lambda + \alpha)z)}{(z^2+\lambda yz +\kappa y^2)^{\gamma}}} \psi_i(y,z)$  & $i=1,2,3,\,\alpha\neq0,$ & $(\lambda \gamma-\mu)X_5 -\mu X_6+\gamma X_7- \kappa \gamma X_8$\\
 \multicolumn{4}{c}%
 {{Here, $\psi_1(y,z)=e^{\displaystyle{\frac{2\lambda \gamma-4\mu}{p}\arctan{\frac{\lambda z+2\kappa y}{pz}}}}$ with $4\kappa-\lambda^2=p^2,\,p\neq0$;}}\\
 \multicolumn{4}{c}%
 {{$\psi_2(y,z)=\displaystyle{\left(\frac{2\kappa y+(\lambda + p)z}{2\kappa y+(\lambda -p)z}\right)^{\frac{2\mu-\lambda\gamma}{p}}}$ with $4\kappa-\lambda^2=-p^2,\,p\neq0$; and}}\\
 \multicolumn{4}{c}%
  {{ $\psi_3(y,z)=e^{\displaystyle{-\frac{4(\mu y+\gamma z)}{\lambda y +2z}}}$ with $4\kappa-\lambda^2=0$
   }}\\

 $f_0 \left(\displaystyle{\frac{y}{y+z}}\right)^{\kappa} y^{1-2\gamma}$ & $ \displaystyle{\left(g_0 - f_0 \frac{y}{y+z}\right)\left(\frac{y}{y+z}\right)^{\kappa-1}} y^{1-2\gamma}$ & $\gamma\neq0,\,\kappa\neq0$ 
    & $\kappa X_2 + 2 (X_6 + X_8)$ \\

 \hline\\
 \hline\\
 \multicolumn{4}{c}%
 {{\bfseries Subalgebra 1. $\gamma X_2 + X_5 + \alpha X_6,\,\alpha\neq0$ }} \\

 $F$ & $G$ & Relations & Additional Extension of Kernel \\

 $f_0 z^{-\kappa}y^{\tilde{\gamma}+1}$ & $g_0 z^{1-\kappa}y^{\tilde{\gamma}}$ &$\tilde{\gamma}=\alpha\kappa-2\gamma,\,\alpha\neq0,\, 1/2,\, 1,\,\kappa\neq0$  
         & $ \kappa X_2 + 2 X_6$ \\
 \hline\\
 \hline\\
 \multicolumn{4}{c}%
 {{\bfseries Subalgebra 2. $\gamma X_2 + X_4 + X_5$ }} \\

 $F$ & $G$ & Relations & Additional Extension of Kernel \\

 $f_0y^{(\kappa+1)}e^{-\alpha z}$ & $g_0y^{(\kappa)}e^{-\alpha z}$ & $\gamma = \displaystyle {\frac{\alpha-\kappa}{2}},\,\kappa\alpha\neq0$ &  $\alpha X_5 +\kappa X_4$ \\ 
 \hline\\
 \hline\\
 \multicolumn{4}{c}%
 {{\bfseries Subalgebra 3. $- X_7 + X_8$}}\\

 $F$ & $G$ & Relations & Additional Extension of Kernel \\
 $(f_0 \cos(u) +g_0 \sin(u))v^{\kappa}$ & $(f_0 \sin(u) -g_0 \cos(u))v^{\kappa}$ & $f_0\neq0,\,g_0\neq0,$ & \\
    &    & $u=\arctan(z/y),\,v^2=y^2+z^2$  &$\displaystyle{\frac{1-\kappa}{2}}X_2 + X_5 + X_6$\\

 \hline\\
 \hline\\
 \multicolumn{4}{c}%
 {{\bfseries Subalgebra 3. $\gamma X_2 - X_7 + X_8,\,\gamma\neq0$}}\\

 $F$ & $G$ & Relations & Additional Extension of Kernel \\
 $e^{\tilde{\gamma}u}(f_0 \cos(u) +g_0 \sin(u))v^{-\tilde{\gamma}\kappa-3}$ & $e^{\tilde{\gamma}u}(f_0 \sin(u) -g_0 \cos(u))v^{-\tilde{\gamma}\kappa-3}$ & $f_0\neq0,\,g_0\neq0,$ & \\
  &  & $\tilde{\gamma}=-2\gamma\neq0,\,u=\arctan(z/y),\,v^2=y^2+z^2$ & $2 X_2 + X_5 + X_6 + \kappa(X_8 - X_7)$\\

 \hline\\
 \hline\\
 \multicolumn{4}{c}%
 {{\bfseries Subalgebra 4. $\gamma X_2 -X_3 +\alpha (X_5+X_6)- X_7 + X_8,\, \alpha>0$}}\\

 $F$ & $G$ & Relations & Additional Extension of Kernel \\
 $e^{(\alpha-2\gamma)u}(f_0 \cos(u) +g_0 \sin(u))v^{\kappa}$ & $e^{(\alpha-2\gamma)u}(f_0 \sin(u) -g_0 \cos(u))v^{\kappa}$ & $f_0\neq0,\,g_0\neq0,$ & \\
    &    & $u=\arctan{\left(\displaystyle{\frac{z+\chi_2(\alpha)}{y-\chi_1(\alpha)}}\right)},$& \\
    &    & $v^2=e^{-2\alpha u}((y-\chi_1(\alpha))^2+(z+\chi_2(\alpha))^2),$  &\\
    &    & $\chi_1(\alpha) = \displaystyle{\frac{\alpha}{\alpha^2+1}},\,\chi_2(\alpha) = \displaystyle{\frac{1}{\alpha^2+1}}$ & $\displaystyle{\frac{1-\kappa}{2}}X_2 + X_5 + X_6 -\chi_1 X_3 + \chi_2 X_4$\\

 \hline\\
 \hline\\
 \multicolumn{4}{c}%
 {{\bfseries Subalgebra 4. $\gamma X_2 +X_3 +\alpha (X_5+X_6)- X_7 + X_8,\, \alpha>0$}}\\

 $F$ & $G$ & Relations & Additional Extension of Kernel \\
 $e^{(\alpha-2\gamma)u}(f_0 \cos(u) +g_0 \sin(u))v^{\kappa}$ & $e^{(\alpha-2\gamma)u}(f_0 \sin(u) -g_0 \cos(u))v^{\kappa}$ & $f_0\neq0,\,g_0\neq0,$ & \\
    &    & $u=\arctan{\left(\displaystyle{\frac{z-\chi_2(\alpha)}{y+\chi_1(\alpha)}}\right)},$& \\
    &    & $v^2=e^{-2\alpha u}((y+\chi_1(\alpha))^2+(z-\chi_2(\alpha))^2),$  &\\
    &    & $\chi_1(\alpha) = \displaystyle{\frac{\alpha}{\alpha^2+1}},\,\chi_2(\alpha) = \displaystyle{\frac{1}{\alpha^2+1}}$ & $\displaystyle{\frac{1-\kappa}{2}}X_2 + X_5 + X_6 +\chi_1 X_3 - \chi_2 X_4$\\

 \hline\\
 \hline\\
 \multicolumn{4}{c}%
 {{\bfseries Subalgebra 4. $\gamma X_2 +\alpha (X_5+X_6)- X_7 + X_8,\, \alpha>0$}}\\

 $F$ & $G$ & Relations & Additional Extension of Kernel \\
 $e^{(\alpha-2\gamma)u}(f_0 \cos(u) +g_0 \sin(u))v^{\kappa}$ & $e^{(\alpha-2\gamma)u}(f_0 \sin(u) -g_0 \cos(u))v^{\kappa}$ & $f_0\neq0,\,g_0\neq0,$ & \\
    &    & $u=\arctan(z/y),\,v^2=e^{-2\alpha u}(y^2+z^2)$  &$\displaystyle{\frac{1-\kappa}{2}}X_2 + X_5 + X_6$\\

 \hline\\
 \hline\\
 \multicolumn{4}{c}%
 {{\bfseries Subalgebra 5. $\gamma X_2 + X_7$ }} \\

 $F$ & $G$ & Relations & Additional Extension of Kernel \\

 $g_0 z^{\beta-1} e^{-y/z} (y+\kappa\tilde{\gamma}z) $ & $g_0 z^{\beta} e^{-y/z} $ & $\tilde{\gamma}=2\gamma$ & 
   $X_5 + \tilde{\gamma} X_6 + (\beta-1)X_7$\\

 \hline\\
 \hline\\
 \multicolumn{4}{c}%
 {{\bfseries Subalgebra 5. $\gamma X_2 + X_4 + X_7$ }} \\

 $F$ & $G$ & Relations & Additional Extension of Kernel \\

 $(g_0 z+f_0)e^{(\beta u-2\gamma z)}$ & $g_0 e^{(\beta u-2\gamma z)}$ & $u=z^2-2y,\, \beta\neq0$ & 
   $ \beta X_2 + X_3  $\\

  \hline\\
 \hline\\
 \multicolumn{4}{c}%
 {{\bfseries Subalgebra 5. $X_4 + X_7$ }} \\

 $F$ & $G$ & Relations & Additional Extension of Kernel \\

 $(g_0 z+ f_0 (\beta + u)^{1/2})(\beta + u)^{\kappa}$ & $g_0 (\beta + u)^{\kappa} $ & $u=z^2-2y,\,\kappa\neq0$ & 
   $\left(1-2\kappa\right) X_2 + 2 \left(- \beta X_3 + 2 X_5 + X_6\right)$\\

 \hline\\
 \hline\\
 \multicolumn{4}{c}%
 {{\bfseries Subalgebra 6. $\gamma X_2 + X_5 + X_6 + X_7$}}\\

 $F$ & $G$ & Relations & Additional Extension of Kernel \\

 $(g_0y+f_0z)z^{\kappa-1}e^{-\tilde{\gamma}(y/z)}$ & $g_0 z^{\kappa}e^{-\tilde{\gamma}(y/z)}$ & $\tilde{\gamma}=2\gamma+k-1\neq0$ &
   $ (\kappa-1)X_2 - 2(X_5 +X_6)$\\ 
 \hline\\
 \hline\\
 \multicolumn{4}{c}%
 {{\bfseries Subalgebra 7. $\gamma X_2 +X_3,\,\gamma\neq0$}}\\

 $F$ & $G$ & Relations & Additional Extension of Kernel \\

 $(f_0 z^{\beta-1} e^{\kappa z-\tilde{\gamma}y})(\kappa z + \tilde{\gamma} \phi_1)$ & $g_0 z^{\beta} e^{\kappa z -\tilde{\gamma} y}$ & $\tilde{\gamma}=2\gamma,\,f_0=g_0/\tilde{\gamma}$&
   $ (\beta-1)X_3 + \kappa X_7 + \tilde{\gamma} X_6$\\
 $(g_0 e^{\beta z + \kappa z^2-\tilde{\gamma}y})\phi(z) $ & $g_0 e^{\beta z + \kappa z^2 -\tilde{\gamma} y}$ & $\tilde{\gamma}=2\gamma,$&  \\
     &   &  $\,\phi=\phi_0z+ \phi_1,\,\phi_0\neq0,\,\kappa=(\tilde{\gamma}\phi_0)/2$& 
   $\beta X_3 + 2\kappa X_7 + \tilde{\gamma} X_4$\\

  \hline
\label{eq:xipp0d}
\end{longtable}
\end{center}
}
\end{landscape}

\begin{appendix}
Appendix 1.
For $\xi^{\prime\prime}=0$, the determining equations have the form
\[
\left\{ \begin{array}{c}
\zeta_{1}^{\prime\prime}=F_{y}\zeta_{1}+F_{z}\zeta_{2}+q_{1},\\
\zeta_{2}^{\prime\prime}=G_{y}\zeta_{1}+G_{z}\zeta_{2}+q_{2},
\end{array}\right.
\]
where $q_{1}$ and $q_{2}$ are functions of $y$ and $z$. Differentiating
them with respect to $x$, one has
\[
\left\{ \begin{array}{c}
\zeta_{1}^{\prime\prime\prime}=F_{y}\zeta_{1}^{\prime}+F_{z}\zeta_{2}^{\prime},\\
\zeta_{2}^{\prime\prime\prime}=G_{y}\zeta_{1}^{\prime}+G_{z}\zeta_{2}^{\prime},
\end{array}\right.
\]
Differentiating the latter equations with respect to $y$ and $z$
\[
\left\{ \begin{array}{c}
F_{yy}\zeta_{1}^{\prime}+F_{yz}\zeta_{2}^{\prime}=0,\\
F_{yz}\zeta_{1}^{\prime}+F_{zz}\zeta_{2}^{\prime}=0,
\end{array}\right.
\]
\[
\left\{ \begin{array}{c}
G_{yy}\zeta_{1}^{\prime}+G_{yz}\zeta_{2}^{\prime}=0,\\
G_{yz}\zeta_{1}^{\prime}+G_{zz}\zeta_{2}^{\prime}=0,
\end{array}\right.
\]
Case 1. Let $F_{zz}\neq0$, then
\[
\zeta_{2}^{\prime}=-\frac{F_{yz}}{F_{zz}}\zeta_{1}^{\prime},\,\,\, F_{yy}-\frac{F_{yz}^{2}}{F_{zz}}=0,\,\,\, G_{yy}-G_{yz}\frac{F_{yz}}{F_{zz}}=0,\,\,\, G_{yz}-G_{zz}\frac{F_{yz}}{F_{zz}}=0.
\]
Thus,
\[
\frac{F_{yz}}{F_{zz}}=k,
\]
and
\[
F_{yy}-kF_{yz}=0,\,\,\, G_{yy}-kG_{yz}=0,\,\,\, G_{yz}-kG_{zz}=0
\]
or
\[
(F_{y}-kF_{z})_{y}=0,\,\,\,(F_{y}-kF_{z})_{z}=0,\,\,\,(G_{y}-kG_{z})_{y}=0,\,\,\,(G_{y}-kG_{z})_{z}=0,
\]
One has
\[
F_{y}-kF_{z}=k_{1},\,\,\, G_{y}-kG_{z}=k_{2},
\]
\[
\frac{dy}{1}=\frac{dz}{-k}=\frac{dF}{k_{1}},\,\,\, F=\Phi(z+ky)+k_{1}y,
\]
\[
\frac{dy}{1}=\frac{dz}{-k}=\frac{dG}{k_{2}},\,\,\, G=\Psi(z+ky)+k_{2}y.
\]
Changing the variables
\[
\bar{y}=y,\,\,\,\bar{z}=z+ky,
\]
the original system
\[
y^{\prime\prime}=F(y,z),\ \ z^{\prime\prime}=G(y,z),
\]
becomes
\[
y^{\prime\prime}=\Phi(\bar{z})+k_{1}y,\,\,\,\bar{z}^{\prime\prime}=\left(k\Phi(\bar{z})+\Psi(\bar{z})\right)+(kk_{1}+k_{2})y
\]
Thus one needs to study the equations
\[
y^{\prime\prime}=k_{1}y+F(z),\,\,\, z^{\prime\prime}=k_{2}y+G(z),\,\,\, k_{2}F^{\prime\prime}\neq0.
\]
\[
\left\{ \begin{array}{c}
\zeta_{1}^{\prime\prime}=k_{1}\zeta_{1}+F^{\prime}\zeta_{2}+q_{1},\\
\zeta_{2}^{\prime\prime}=k_{2}\zeta_{1}+G^{\prime}\zeta_{2}+q_{2},
\end{array}\right.
\]

Notice that because of $F^{\prime\prime}\neq0$, then $\zeta_{2}^{\prime}=0$
from the second equation
\[
0=k_{2}\zeta_{1}^{\prime}\Rightarrow\zeta_{1}^{\prime}=0.
\]

Case 2. Let $F_{zz}=0$, then by symmetry $G_{yy}=0$. Hence
\[
F_{yy}\zeta_{1}^{\prime}+F_{yz}\zeta_{2}^{\prime}=0,\,\,\, F_{yz}\zeta_{1}^{\prime}=0,\,\,\, G_{yz}\zeta_{2}^{\prime}=0,\,\,\, G_{yz}\zeta_{1}^{\prime}+G_{zz}\zeta_{2}^{\prime}=0,
\]
If
\[
F_{yz}\neq0\Rightarrow\zeta_{1}^{\prime}=0,\,\,\,\zeta_{2}^{\prime}=0.
\]
Hence,
\[
F_{yz}=0,\,\,\, F_{zz}=0,\,\,\, G_{yy}=0,\,\,\, G_{yz}=0
\]
and
\[
F_{yy}\zeta_{1}^{\prime}=0,\,\,\, G_{zz}\zeta_{2}^{\prime}=0.
\]
Thus
\[
y^{\prime\prime}=k_{1}z+F(y),\,\,\, z^{\prime\prime}=k_{2}y+G(z),\,\,\, k_{1}k_{2}(F^{\prime\prime2}+G^{\prime\prime2})\neq0.
\]
\[
\left\{ \begin{array}{c}
\zeta_{1}^{\prime\prime}=F^{\prime}\zeta_{1}+k_{1}\zeta_{2}+q_{1},\\
\zeta_{2}^{\prime\prime}=k_{2}\zeta_{1}+G^{\prime}\zeta_{2}+q_{2},
\end{array}\right.
\]

Let $F^{\prime\prime}\neq0$, then $\zeta_{1}^{\prime}=0$ and because
\[
\zeta_{1}^{\prime\prime\prime}=F^{\prime}\zeta_{1}^{\prime}+k_{1}\zeta_{2}^{\prime}\Rightarrow k_{1}\zeta_{2}^{\prime}=0\Rightarrow\zeta_{2}^{\prime}=0
\]

Let $G^{\prime\prime}\neq0$, then $\zeta_{2}^{\prime}=0$ and because
\[
\zeta_{2}^{\prime\prime\prime}=k_{2}\zeta_{1}^{\prime}+G^{\prime}\zeta_{2}^{\prime}\Rightarrow k_{2}\zeta_{1}^{\prime}=0\Rightarrow\zeta_{1}^{\prime}=0
\]
\end{appendix}
\end{document}